\documentclass[11pt]{amsart}
\usepackage{cite}
\usepackage{amsbsy,amsfonts,amsmath,amssymb,amscd,amsthm}
\usepackage{mathrsfs}
\usepackage{subfigure,enumitem,graphicx}
\usepackage[a4paper,text={6.1in,8in},centering,includefoot,foot=1in]{geometry}
\usepackage[colorlinks=true,linkcolor=black,citecolor=black]{hyperref}
\usepackage{graphics}
\usepackage{tikz}
\usetikzlibrary[patterns]
\def\leb{\ensuremath{\text{Leb}}}
\def\sign{\ensuremath{\text{sgn}}}

\small\normalsize
\theoremstyle{plain}
\newtheorem{theorem}{Theorem}
\newtheorem{maintheorem}{Theorem}

\newtheorem{lemma}[theorem]{Lemma}

\newtheorem*{theorem*}{Theorem}

\newtheorem{example}[theorem]{Example}
\newtheorem{proposition}[theorem]{Proposition}

\newtheorem*{claim*}{Claim}

\newtheorem*{conj*}{Conjecture}

\let\oldtocsection=\tocsection
\let\oldtocsubsection=\tocsubsection
\renewcommand{\tocsection}[2]{\hspace{0em}\bf\oldtocsection{#1}{#2}}
\renewcommand{\tocsubsection}[2]{\hspace{1em}\oldtocsubsection{#1}{#2}}
\def\leb{\ensuremath{\text{Leb}}}
\def\sign{\ensuremath{\text{sgn}}}
\def\diff{\textsf{Diff}}

\begin{document}
\title{Saturation of Generalized Partially Hyperbolic Attractors}
\author{Abbas Fakhari}
\address{Department of Mathematics, Shahid Beheshti University,
Tehran 19839, Iran\,\,\,\&\newline
School of Mathematics, IPM,
Tehran 19395-5746, Iran.} \email{a\_fakhari@sbu.ac.ir}
\subjclass[2000]{37DXX, 37XX}
\author{Mohammad Soufi}
\address{Instituto de Matem\'atica e Estat\'istica, Universidade do Estado do Rio de Janeiro, Rua S\~ao Francisco Xavier, 524, Pavilh\~ao Reitor Jo\~ao Lyra Filho - 6 andar, Bloco B, 20550-900,
Rio de Janeiro-RJ, Brasil } \email{msoufin@gmail.com} \keywords{generalized partially hyperbolic attractor, $C^1$-generic, horseshoe-like, Lorenz attractor}

\begin{abstract}
We prove the saturation of a generalized partially hyperbolic attractor of a  $C^2$ map. As a  consequence, we show that any generalized partially hyperbolic horseshoe-like attractor of a $C^1$-generic diffeomorphism has zero volume. In contrast, by modification of the Poincar\'e cross section of Lorenz geometric model, we build a $C^1$-diffeomorphism with a partially hyperbolic horseshoe-like attractor of positive volume.
\end{abstract}

\maketitle

\section{introduction}
A hyperbolic set is a compact invariant set over which the tangent bundle splits into two invariant sub-bundles, one is contracting and the other one is expanding. The Lebesgue measure (volume) of hyperbolic sets is an interesting subject considered in many articles. The scenario begun by the seminal works of Bowen in 70's. Bowen proved in \cite{B2} that a hyperbolic invariant set of positive volume of a $C^2$-diffeomorphism does contain some stable and unstable manifolds. On the other hand, he showed in \cite{B1} the existence of a $C^1$-diffeomorphism admitting a totally disconnected hyperbolic set of positive volume. The issue of the volume and the interior of a hyperbolic set followed by many
authors. For instance, it is shown in \cite{AP} that a transitive hyperbolic set which attracts a set with positive volume necessarily attracts a neighborhood of itself. It is also proved in \cite{AP} that there are no proper transitive hyperbolic sets with positive volume for diffeomorphisms whose differentiability is higher than one. In the context of volume preserving diffeomorphisms, it is proved in \cite{BV} that a volume preserving diffeomorphism with a
hyperbolic set of positive volume should be an Anosov diffeomorphism (see also \cite{X}). The main point in this context is the saturation of an invariant set of positive volume with stable or unstable leaves. In light of the
rich consequences on the volume of a hyperbolic set, a lot of studies has been undertaken in a more general landscape which are non-uniform hyperboliciy, partial hyperbolicity and maps with some discontinuities. For instance, Alves and Pinheiro proved in \cite{AP} the non-existence of
horseshoe-like partially hyperbolic sets with positive volume for $C^{2}$-diffeomorphims. Extending their result, Zhang  proved in \cite{Z} that the hyperbolic lamination over the limit sets of s-density points of a partially hyperbolic set with positive volume of a $C^2$-diffeomorphism should be contained in the partially hyperbolic set. The aim of this article is to study
the volume of invariant sets in the context of  non-uniform and partially hyperbolic dynamics with some discontinuities. First, we prove that a generalized hyperbolic attractor of a $C^2$-diffeomorphism with positive volume should contain the stable and unstable laminations. Then inspired by the Bowen's example of a fat horseshoe, we build a $C^1$-flow whose time one map admits a horseshoe-like generalized partially hyperbolic attractor with positive volume. Our construction is based on a modification of a classical Lorenz flow using the Bowen map. We begin by recalling some essential notions.


\section{Generalized Partially Hyperbolic Attractors}
Let $M$ be a smooth Riemannian manifold, $U$ an open bounded connected subset of $M$ with compact closure and $\Sigma$ a closed subset of $U$. 
Let $f:U\setminus \Sigma\to U$ be a twice differentiable map satisfying conditions (H1) and (H2);
\begin{itemize}
\item[\bf{H1)}] $f$ is a $C^2$-diffeomorphism from the open set $U \setminus \Sigma$ onto its image $f(U \setminus \Sigma)$.
\item[\bf{H2)}] there exists $C_i > 0$ and $\alpha_i \geq 0 ~(i=1,2)$ such that
\begin{align*}
\left\| d^2 f_x \right\| &\leq C_1 d(x,\Sigma^+)^{-\alpha_1} \quad \text{for any } x\in U \setminus \Sigma,\\
\left\| d^2 f_x^{-1} \right\| &\leq C_2 d(x,\Sigma^-)^{-\alpha_2} \quad \text{for any } x\in f(U \setminus \Sigma),
\end{align*}
where $d$ is the Riemannian distance in $M$, $\Sigma^+=\Sigma\cup
\partial(U)$ and $\Sigma^-=\partial(f(U\setminus \Sigma))$.
\end{itemize}
Let $U^+=\{z\in U: f^n(z)\not\in \Sigma^+\text{~for all~} n\geq 0\}$, $D=\bigcap_{n\geq 0} f^n(U^+)$ and $\Lambda=\overline{D}$. The set $\Lambda$ is called an {\it attractor} for $f$. 
For $\epsilon > 0$ and $l \in \mathbb{N}$, put
\begin{align*}
D^+_{\epsilon,l}&=\{z\in \Lambda : d(f^{n}(z),\Sigma^+)\geq l^{-1}e^{-\epsilon n},~n=0,1,\ldots\},\\
D^-_{\epsilon,l}&=\{z\in \Lambda : d(f^{-n}(z),\Sigma^-)\geq l^{-1}e^{-\epsilon n},~n=0,1,\ldots\},\\
D^0_{\epsilon,l}&=D^+_{\epsilon,l}\cap D^-_{\epsilon,l}\,\,\,\,\text{and}\,\,\,\, D^0_{\epsilon}=\bigcup_{l\geq 0} D^0_{\epsilon,l}.
\end{align*}
Suppose that for any $z\in U$, there are two cones
$C^s_\alpha(z)$ and $C^c_\alpha(z)$ such that  for some
$0<\lambda<1$ and $\lambda<\mu$,
\begin{enumerate}
\item the angel between $C^s_\alpha(z)$ and $C^c_\alpha(z)$ is uniformly bounded away from zero,
\item \begin{enumerate}
      \item $df(C^c_\alpha(z))\subset C^c_\alpha(f(z))$, for $z\in U\setminus\Sigma^+$,
      \item $df^{-1}(C^s_\alpha(z))\subset C^s_\alpha(f^{-1}(z))$, for $z\in f(U\setminus \Sigma^+)$,
      \end{enumerate}
\item 
\begin{enumerate}
      \item  $\|df^n(z)v\|\geq C \mu^n\|v\|$ for any $v\in C^c_\alpha(z)$ and $z\in U^+$,
      \item  $\|df^{-n}(z)v\|\geq C \lambda^{-n}\|v\|$ for any $v\in C^s_\alpha(z)$ and $z\in f^n(U^+)$.
      \end{enumerate}
\end{enumerate}
For any $z\in D$, put
$$E^s(z)=\bigcap_{n\geq 0} df^{-n}\big( C^s_\alpha(f^n(z))\big)\quad\text{and}\quad
E^c(z)=\bigcap_{n\geq 0} df^{n}\big( C^c_\alpha(f^{-n}(z))\big).$$
Subspaces $E^s$ and $E^c$, which are complement to each other, are invariant and
for any $z\in \bigcap_{n\geq 0} f^n(U^+)$ and for any $w\in E^s(z)$ and $v\in E^c(z)$,
$$\|df^n(z)w\|\leq C\lambda^n\|w\|\quad\text{and}\quad \|df^n(z)v\|\geq C\mu^n\|v\|.$$
\begin{proposition}\label{pro:lsm}
There is $\epsilon>0$ such that for any $z\in
D^+_{\epsilon,l}$ there is a $C^1$-submanifold $W^s_{\delta_l}(z)$
of size $\delta_l$, called local stable manifold, such that
\begin{enumerate}
\item $T_zW^s_{\delta_l}(z)=E^s(z)$,
\item $f(W^s_{\delta_l}(z))\subset W^s_{\delta_{l_1}}(f(z))$, where $f(z)\in D^+_{\epsilon,l_1}$,
\item there is $\lambda<\tilde\lambda<1$ and $C_l>0$ such that for $x,y\in W^s_{\delta_l}(z)$, $$d(f^n(x),f^n(y))\leq C_l\tilde\lambda^n d(x,y),$$
\item $W^s_{\delta_l}(z)$ varies continuously on $D^+_{\epsilon,l}$~,
\item there is some  $\gamma > 0$ such that $\delta_{l(f^m(z))}\geq \exp(-\gamma|m|)\delta_l$ for $m\in\mathbb{Z}$.
\end{enumerate}
\end{proposition}
\begin{proof}
See Proposition 4 in \cite{p}.
\end{proof}
An invariant subset of $\Lambda$ is {\it horseshoe-like} if it doesn't contain any connected piece of
stable manifold. In general, the attractor $\Lambda$ may has empty interior or even has zero Lebesgue measure.
Here, we focus on the case in which the attractor has positive Lebesgue measure.
\begin{maintheorem}
\label{thma1} If $D_\epsilon^0$ supports an ergodic absolutely
continuous invariant measure $\mu$ then for any invariant set of
positive Lebesgue measure $A$ in $D_{\epsilon}^0$ we have
\begin{center}
$W^s_{loc}(x)\subset \overline{A}\quad$ for $\mu$-a.e $x\in A$.
\end{center}
\end{maintheorem}
The following example realizes the condition of Theorem \ref{thma1}.
\begin{example}
Consider the map $F$ on the unite square defined by
$$(x,y)\mapsto (2x, \lambda y+f(x))$$
where $f$ is a positive $C^r$ $(r\geq 1)$ function and $\lambda<1$. The function $F$ is an Anosov endomorphism with the singular set $\Sigma=\{(x,y):x=1/2\}$. The mapping $F$ admits an
ergodic SRB measure $\mu$ with the full basin supported on an attractor $D$. In the case when $\lambda<1/2$, The attractor may have zero Lebesgue measure and hence the SRB measure is totally singular with respect to the Lebesgue measure. The case when $\lambda\geq1/2$ is more interesting as it admits, in a generic way, an ergodic SRB measure $\mu$ which is absolutely continuous with respect to the Lebesgue measure. Avila et al. proved in \cite{agt} that for suitable choice of $r$ and $\lambda$, the density function $g=\frac{d\mu}{dx}$ is essentially bounded and the attractor $D$ has non-empty interior. In particular one has that
\begin{equation*}
\left|\int_D \log d(z,\Sigma^+)~d\mu(z)\right|\leq \int_D |\log d(z,\Sigma^+)g(z)|~dz\leq \|g\|_{ess} \int_{[0,1]^2} |\log d(z,\Sigma^+)|~dz<\infty,
\end{equation*}
which implies the existence of a constant $\varepsilon_0$ that $\mu(D^0_\epsilon)>\varepsilon_0$ for any $\epsilon>0$, see \cite[Proposition 2]{p}. Hence, by the absolute continuity, we have $\leb(D_{\epsilon}^0)>0$ for any $\epsilon$.
\end{example}
%
%

Next goal is to complement our scenario by presenting a $C^1$ map with
a generalized hyperbolic set of positive Lebesgue measure which doesn't contain
any stable leaf. 
%
The theorem \ref{thmb} gives an example similar to the Bowen's example of horseshoe with positive volume
in the context of generalized hyperbolic attractors. The existence of this example has been claimed in [1].
Here we give a complete proof using a different approach which is  based on a modification of classical geometric
model for Lorenz attractor introduced in \cite{gh} for the first time.
\begin{maintheorem}
\label{thmb} There exist a $C^1$-flow on a three
dimensional manifold with a singular hyperbolic attractor, such that its time one map is a $C^1$-diffeomorphism admitting a horseshoe-like partially hyperbolic
attractor of positive volume.
\end{maintheorem}
\section{Hyperbolicity and Saturation}
In this section we are going to prove our theorems concerning the
saturation of the invariant set of a $C^2$-diffeomorphism having weaker
kind of hyperbolicity discussed above. We begin by an essential
ingredient. In the rest of the paper we assume that the attractor is regular,
i.e. for all small $\epsilon$, the set $D^0_{\epsilon}$ is non-empty, see \cite[condition H3)]{p}.
In \cite{X}, the author introduced a simple dynamical density basis
for positive Lebesgue measure sets, which is strongly inspired by
\cite{bw}, called s-density point. 
 Let $A$ be a measurable subset of $D^0_\epsilon$. A point $z\in A\cap D^0_{\epsilon,l}$
is called a $s${\it -density point} of $A$ if
$$\lim_{n\to\infty} \frac{m_s(B_n^s(z)\cap A)}{m_s(B_n^s(z))}=1,$$
where $m_s$ is Lebesgue measure on $W^s_{\delta_l}(z)$ and $B_n^s(z)= f^n\left(W^s_{\delta_{l(n)}}\left(f^{-n}(z)\right)\right)$. 
To guarantee the existence of such s-density points for sets of
positive measure, we need to examine two basic properties of the
shrinking sequence $B_n^s(z)$: scaling and engulfing.

\hspace{-.4cm}{\bf Scaling:} For any integer $k$, there is a
constant $C=C(k)$ such that for 
all $z\in D_{\epsilon,l}^0$,
$$
\frac{m_s(B^s_{n+k}(z))}{m_s(B^s_n(z))}\geq C.
$$
Let us first prove the bounded distortion property which is the key ingredient for proving the scaling property.
\begin{lemma}[Bounded Distortion] \label{prop-dist-2}
There exists a constant $L_{\mathcal{H}}=L_{\mathcal{H}}(l)$ such that for any $z\in
D^0_{\epsilon,l}$ and any $x,y\in W^s_{\delta_l}(z)$,
$$
L_{\mathcal{H}}<\bigg|\frac{J f^n(x)}{J f^n
(y)}\bigg|<L_{\mathcal{H}}^{-1},\quad\text{for all}~n\in\mathbb{N}.
$$
\end{lemma}
\begin{proof}
Note that by assumption, $\digamma=\log |J f|$ is $\alpha$-H\"older
and hence, for every $x,y$, it holds $|\digamma(x)-\digamma(y)|\leq
C d(x,y)^\alpha$ for some constant $C>0$. Let $z\in
D^0_{\epsilon,l}$ and $x,y\in W^s_{\delta_l}(z)$. 
Then it holds,
$ d(f^i(x),f^i(y))\leq C_l \tilde\lambda^id(x,y)$. This implies that
\begin{align*}
\log\frac{|J f^n(x)|}{|J f^n(y)|}
\leq\sum_{i=0}^{n-1}|\digamma(f^i(x))-\digamma(f^i(y))| \leq C
\sum_{i=0}^{n-1} d(f^i(x),f^i(y))^\alpha &\leq C
\sum_{i=0}^{n-1}(C_l\tilde\lambda^i )^\alpha\\ &\leq C
C_l^\alpha\sum_{i=0}^{\infty}\tilde\lambda^{i\alpha}.
\end{align*}
Now, taking $L_{\mathcal{H}}=\exp\{CC_l^\alpha/ (1-\tilde\lambda^{\alpha}) \}$, the desired inequalities hold.
\end{proof}
\begin{proof}[Proof of Scaling]
Suppose that  $z\in D_{\epsilon,l}^0$. We have
\begin{align*}
m_s\left(B_n^s(z)\right)=m_s\left(f^n\big(W^s_{\delta_{l(n)}}(f^{-n}(z))\big)\right)=&\int_{W^s_{\delta_{l(n)}}(f^{-n}(z))}
|Jf^n(y)|~dm_s(y)\\
\leq& \max_{y\in W^s_{\delta_{l(n)}}(f^{-n}(z))}
|Jf^n(y)|~\delta_{l(n)}.
\end{align*}
With the same argument, it holds
\begin{align*}
m_s(B_{n+k}^s(z))=&m_s\left(f^{n+k}\big(W^s_{\delta_{l(n+k)}}(f^{-n-k}(z))\big)\right)\geq
\min_{y\in W^s_{\delta_{l(n+k)}}(f^{-n-k}(z))} |J
f^{n+k}(y)|~\delta_{l(n+k)}\\
&\hspace{3.8cm}\geq
\min_{x\in U\setminus\Sigma} |Jf^k(y)|~\min_{y\in W^s_{\delta_{l(n)}}(f^{-n}(z))} |J
f^{n}(y)| ~\delta_{l(n+k)}.
\end{align*}
By the fifth item of Proposition \ref{pro:lsm}, there is a constant
$C_2=C_2(k)$ such that $\delta_{l(n+k)}\geq C_2\delta_{l(n)}$.
On the other hand, by the bounded distortion property, one has
$$
L_{\mathcal{H}}<\frac{\max_{y\in
W^s_{\delta_{l(n)}}(f^{-n}(z))} |J f^n(y)|}{\min_{y\in
W^s_{\delta_{l(n)}}(f^{-n}(z))} |J
f^{n}(y)|}<L_{\mathcal{H}}^{-1}.
$$
Assembling the two facts above, one gets
$$
\frac{m_s(B^s_{n+k}(z))}{m_s(B^s_n(z))}\geq
\frac{L_{\mathcal{H}}}{C_2\min_{z\in U\setminus\Sigma} |Jf^k(z)|}.
$$
\end{proof}
\hspace{-.4cm}{\bf Engulfing:} For any $l$, there is $L=L(l)$ such that for any
$n\geq 1$ and for any $x,y\in D_{\epsilon,l}^0$, if
$B^s_{n+L}(x)\cap B^s_{n+L}(y)\neq\emptyset$, then it holds
$B^s_{n+L}(x)\cup B^s_{n+L}(y)\subset B^s_n(x)$.
\begin{proof} We show that for any $l$ there is a uniform $L$ such that for any $x,y\in  D_{\epsilon,l}^0$, $$f^L\left(W^s_{\delta_{l(n+L)}}\left(f^{-n-L}(x)\right) \cup W^s_{\delta_{l(n+L)}}\left(f^{-n-L}(y)\right)\right) \subset W^s_{\delta_{l(n)}}(f^{-n}(x)).$$ For this, it is enough to note that
\[\textrm{Diam}(
f^L \left(
W^s_{\delta_{l(n+L)}}\left(f^{-n-L}(x)\right) \cup W^s_{\delta_{l(n+L)}}\left(f^{-n-L}(y)\right)
\right)
)
\leq C_3C_l\tilde{\lambda}^L\textrm{Diam}\big(W^s_{\delta_{l}}(f^{-n}(x))\big),\]
%
for some constant $C_3$.
\end{proof}
Now, we can state the following folklore proposition on the
measure theoretical behavior of the stable manifolds of s-density
points in our context.
\begin{theorem}\label{mthm}
If $A\subset D^0_\epsilon$ has positive volume then, for any
$x\in A\cap D_{\epsilon,l}^0$ and $\eta>0$, there is $n_0$
depending only to $x$ and $\eta$
, and not to $\delta_l$, such that for any $n\geq n_0$ with 
$$m_s(W^s_{\delta_{l(n)}}(f^{-n}(x)\cap A))\geq
(1-\eta)m_s(W^s_{\delta_{l(n)}}(f^{-n}(x))$$
\end{theorem}
\begin{proof}
With the same calculation above, one has
\[\frac{m_s(B^s_n(x)\setminus A)}{m_s(B^s_n(x))}\geq \frac{\min_{y\in
W^s_{\delta_{l(n)}}(f^{-n}(x))}Jf^n(y)}{\max_{y\in
W^s_{\delta_{l(n)}}(f^{-n}(x))}Jf^n(y)}\cdot\frac{m_s(W^s_{\delta_{l(n)}}(f^{-n}(x))\setminus
A)}{m_s(W^s_{\delta_{l(n)}}(f^{-n}(x)))}.\] The left hand
side of the inequality tends to zero, by the assumption, and so
the right is the same.
\end{proof}


\begin{proof}[proof of Theorem \ref{thma1}]
Suppose that $D^0_\epsilon$ supports an ergodic
absolutely continuous invariant measure $\mu$ and $A\subset
D^0_\epsilon$ is invariant set of positive Lebesgue measure.  By
the ergodicity, $\mu(D_\epsilon^1)=\mu(A)=1$. 
If $A^d$ be the set of s-density points of $A$ then
$Leb(A\setminus A^d)=0$ and hence $\mu(A\setminus A^d)=0$, see \cite{X}.
Thus we can assume that any $\mu$-generic point of $A$ is
s-density point.
Now, let $x\in A$ be a generic point. Choose a sequence $\tau^{(n)}$ of return times of $x$ to $D^0_{\epsilon,l}$
such that $f^{-\tau^{(n)}}(x)\to x$.
By the
continuity of stable manifolds on $D^0_{\epsilon,l}$ , $$W^s_{\delta_{l}}(f^{-\tau^{(n)}}(x))\to
W^s_{\delta_l}(x).$$
 Now, by contradiction, suppose that
$W^s_{\delta_l}(x)\not\subset \overline{A}$. Hence, there is a
point $z\in W^s_{\delta_l}(x)\setminus \overline{A}$. So, for some
$\delta_0\leq\delta_l$, $B_{\delta_0}(z)\cap A=\emptyset$. Let
$x_n\in W^s_{\delta_{l}}(f^{\tau^{(n)}}(x))$ such that $x_n\to z$. For any
$n$, $x_n\in B_{\delta_0/2}(z)$ and so $B_{\delta_0/2}(x_n)\cap
A=\emptyset$, in particular, $W^s_{\delta_0/2}(x_n)\cap
A=\emptyset$. Since $W^s_{\delta_0/2}(x_n)\subset
W^s_{\delta_{l}}(f^{\tau^{(n)}}(x))$, in view of Theorem \ref{mthm}, for any
$\eta>0$ and large enough $\tau^{(n)}$, one gets
\begin{align*}
(1-\eta)\mu^s(W^s_{\delta_{l}}(f^{-\tau^{(n)}}(x)))&\leq
\mu^s(W^s_{\delta_{l}}(f^{-\tau^{(n)}}(x))\cap A)\\&\leq
\mu^s(W^s_{\delta_{l}}(f^{-\tau^{(n)}}(x))\setminus
W^s_{\delta_0/2}(x_n))\\ &=
\mu^s(W^s_{\delta_{l}}(f^{-\tau^{(n)}}(x)))-C_4\delta_0,
\end{align*}
for some constant $C_4$. This is a contradiction since $\tau^{(n)}$ does not depend on $\eta$.
\end{proof}
\section{horseshoe-like generalized partially hyperbolic attractor of positive volume}
\begin{figure}[ht]
\def\svgwidth{8cm}
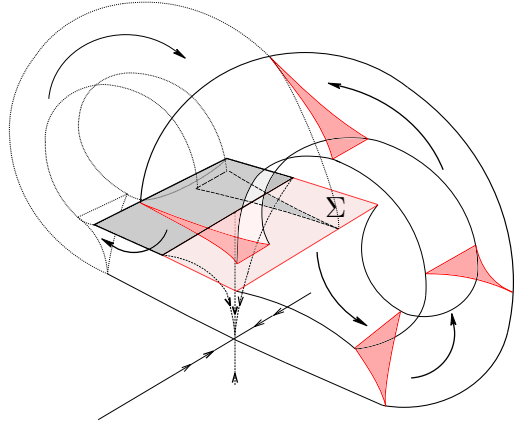
\caption{Geometric model of Lorenz attractor and its Poincar\'e section }
\label{Lorenz}
\end{figure}

In this section we provide an example of generalized hyperbolic set
which is a two dimensional reduction of a singular hyperbolic
attractor. 
Our aim is to present an example showing that the $C^2$ regularity is needed in Theorem \ref{thma1}. The first example of
singular hyperbolic attractor was given by E. Lorenz in 70's during
his study on the foundations of long range weather forecast. In
fact, he gave an ODE that the numerical experiments showed
the corresponding flow exhibits an attractor. Albeit the simplicity
of the Lorenz equations (2-degree polynomial), it was not a
simple task to solve. There are two conceptual problems. The first
one is the presence of a hyperbolic equilibrium accumulated by
regular orbits which prevents the Lorenz attractor from being
hyperbolic and the second one is that the solutions slow down when passing near the equilibrium, which means unbounded return times
and thus unbounded integration errors. The impossibility of solving
the equations leads Afraimovich-Bykov-Shilnikov and
Guckenheimer-Williams, independently (in the seventies), to propose
a geometrical model. In Figure~\ref{Lorenz}, we illustrate schematic
geometric model, Poincar\'e section $\Sigma$ showed by a square
foliated by contracting invariant leaves. The Poincar\'e section
will partition by the intersection of the two-dimensional stable
manifold of the equilibrium in two rectangles whose images after
coming back to the cross section are two triangles  with curved
sides as shown in the Figure \ref{Lorenz}. This model has a transitive
attractor $\Lambda$ with a dense set of periodic points.
\vspace{-9.5cm}
In what follows, $\leb_n$ stands for the Lebesgue measure on $\mathbb{R}^n$ and $\leb=\leb_1$. In this section an example of singular hyperbolic
attractor with positive volume is given. The idea is to construct a fat
horseshoe\footnote{We call a horseshoe $H$ fat if $\leb_2(H)>0$.
}~$H$ on the
Poincar\'e section $\Sigma$ of a geometric Lorenz attractor
$\Lambda$ such that $H\subset\Lambda \cap \Sigma$. Let us assume
that we have constructed the fat horseshoe $H$. Then
$$
H_{\delta}=\left\{X^t(x): |t|\leq\delta, x\in H\right\},
$$
is a subset of $\Lambda$ and
$\leb_3(H_{\delta})\approx\delta\times\leb_2(H)>0$. Therefore the
geometric Lorenz attractor $\Lambda$ has positive volume for its
associated vector filed and a priori for the time one map
diffeomorphism. Note that the attractor $\Lambda$ is horseshoe-like,
i.e., it dose not contain any local stable manifold, but it dose
contain an unstable manifold\footnote{The attractor contains the unstable manifold of the hyperbolic equilibrium point.}.
\subsection{Fat triangles on $\Sigma$}
First, by a straightforward
calculation, it is shown that fat triangles on the cross section
does not guarantee the positiveness of the volume of the cross sectional attractor $\Lambda \cap \Sigma$ which consists of two Cantor cones.

\begin{figure}[h]
\centering
 \begin{tikzpicture}
    \draw (0,0) rectangle (4,4);
  \draw[thick] (2,0) -- (2,4);
  \draw (3,0) -- (3,4);
  \node at (3,-.3) {$x=a$};
  \node at (2,-.3) {$\Gamma$};
  \draw (4,3.5) .. controls (2,3.2)  .. (0,3);
    \draw (4,3.4) .. controls (2,3.16) .. (0,3);
      \draw (4,3.3) .. controls (2,3.12) .. (0,3);
        \draw (4,3.2) .. controls (2,3.08) .. (0,3);
          \draw (4,3.1) .. controls (2,3.04) .. (0,3);
            \draw (4,3)   .. controls (2,3)    .. (0,3);
              \draw (4,2.9) .. controls (2,2.96) .. (0,3);
                \draw (4,2.8) .. controls (2,2.92) .. (0,3);
                  \draw (4,2.7) .. controls (2,2.88) .. (0,3);
                    \draw (4,2.6) .. controls (2,2.84) .. (0,3);
                      \draw (4,2.5) .. controls (2,2.8)  .. (0,3);

   \draw (0,1.5) .. controls (2,1.2) .. (4,1);
    \draw (0,1.4) .. controls (2,1.16) .. (4,1);
      \draw (0,1.3) .. controls (2,1.12) .. (4,1);
        \draw (0,1.2) .. controls (2,1.08) .. (4,1);
          \draw (0,1.1) .. controls (2,1.04) .. (4,1);
            \draw (0,1)   .. controls (2,1)    .. (4,1);
              \draw (0,.9) .. controls (2,.96) .. (4,1);
                \draw (0,.8) .. controls (2,.92) .. (4,1);
                  \draw (0,.7) .. controls (2,.88) .. (4,1);
                    \draw (0,.6) .. controls (2,.84) .. (4,1);
                      \draw (0,.5) .. controls (2,.8)  .. (4,1);
      (5,0) .. controls (6,0) and (6,1) .. (5,2);
\end{tikzpicture}
\caption{Sectional attractor (Cantor cones)}
\label{fig:3}
\end{figure}
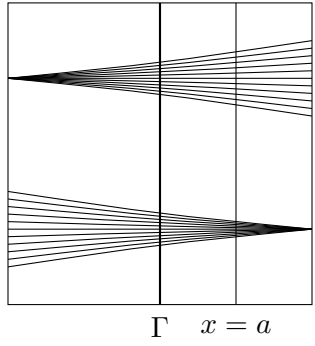

Let $\Sigma=[-1,1]\times [-1,1]$ and
$\Gamma=\{(x,y)\in\Sigma:x=0\}$. For $k\geq 2$, define the function
$F_k(x,y)$ on $\Sigma\setminus \Gamma$ as follows:
\begin{equation*}
 F_k(x,y)=\left\{
\begin{array}{ll}
\left(2x^{1/2}-1,1/2(yx^{1/k}+1)\right)&\quad x>0,\\\\
\left(-2|x|^{1/2}+1,1/2(y|x|^{1/k}-1)\right)&\quad x<0.
\end{array}
\right.
\end{equation*}
We try to calculate the Lebesgue measure of the Cantor set $C(a)$ which
is the intersection of the sectional attractor
$\Lambda\cap\Sigma=\bigcap\limits_{n\geq 0}F_k^n(\Sigma\setminus\Gamma)$ with segment
$x=a$, see Figure \ref{fig:3}.
It can be seen that
\[
C(a)=\lim\limits_{n\to +\infty}\bigcup\limits_{x_0\in f^{-n}(a)}F_k^n(x=x_0),\,\,
\textrm{where},\,
f(x)=\left\{
\begin{array}{ll}
2x^{1/2}-1&\quad x>0,\\\\
-2|x|^{1/2}+1&\quad x<0.
\end{array}
\right.
\]
Let $a_{0,1}=a$ and $b_0=1$. A simple calculation shows that
$$f^{-n}(a)=\left\{\frac{a_{n,1}}{b_n},\frac{a_{n,2}}{b_n},\frac{a_{n,3}}{b_n},\ldots,\frac{a_{n,2^n}}{b_n}\right\},$$
where $a_{n,m}=(-1)^m(a_{n-1,\lceil m/2 \rceil}+(-1)^m b_{n-1})^2$
for $m=1,2,3,\ldots, 2^n,$ and $b_n=2^{2^{n+1}-2}=4b_{n-1}^2$. Using
the notations $A_{n,m}=F_k^n(x=\frac{a_{n,m}}{b_n})$ and $C_n(a)=\bigcup\limits_{m=1}^{2^n}A_{n,m}$, we have

\begin{align*}
\leb(C_n(a))&=\sum\limits_{m=1}^{2^n}\leb(A_{n,m})\\
&=\sum_{\substack{1\leq m\leq 2^n\\m \text{~even}}}\leb(A_{n,m-1})+\leb(A_{n,m})\\
&=\sum_{\substack{1\leq m\leq 2^n\\m \text{~even}}} \frac{1}{2}\left(\left|\frac{a_{n,m-1}}{b_n}\right|^{1/k}+\left|\frac{a_{n,m}}{b_n}\right|^{1/k}\right)\leb(A_{n-1,m/2})\\
&=\sum_{\substack{1\leq m\leq 2^n\\m \text{~even}}}
\frac{1}{2}\left(\frac{(b_{n-1}-a_{n-1,m/2})^{2/k}}{4^{1/k}b_{n-1}^{2/k}}+
\frac{(b_{n-1}+a_{n-1,m/2})^{2/k}}{4^{1/k}b_{n-1}^{2/k}}\right)\leb(A_{n-1,m/2})\\
&=\sum\limits_{m=1}^{2^{n-1}} \frac{1}{2^{1+2/k}}\left(\left(1-\frac{ a_{n-1,m}}{b_{n-1}}\right)^{2/k}+\left(1+\frac{a_{n-1,m}}{b_{n-1}}\right)^{2/k}\right)\leb(A_{n-1,m})\\
&\leq \sum\limits_{m=1}^{2^{n-1}} \frac{1}{2^{2/k}}\leb(A_{n-1,m})\\
&=\frac{1}{2^{2/k}}\leb(C_{n-1}).
\end{align*}
Therefore $\leb(C_n(a))\leq
\frac{1}{4^{n/k}}\leb(C_0(a))=\frac{2}{4^{n/k}}$. Hence, $\leb(C(a))=0$. Then, Fubini theorem implies $\leb_2(\Lambda\cap\Sigma)=0$.
\subsection{Modification of Lorenz map using Bowen's method}
Now, we are trying to give an example of singular hyperbolic
attractor with positive volume. We recall that $f:[-1,1]\to [-1,1]$
is the {\it Lorenz map} if it satisfies:
\begin{itemize}
 \item $f(-1)>-1$, $f(1)<1$ and $\lim\limits_{x\to 0^\pm}f(x)=\mp 1$,
 \item $f'(x)>\alpha>0$ for $x\not= 0$ and $\lim\limits_{x\to 0^\pm}f'(x)=+\infty$.
\end{itemize}

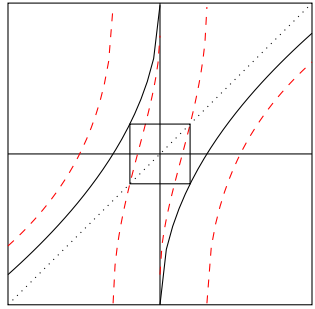
\begin{figure}[h]
\centering
\begin{tikzpicture}[scale=1]
\draw (-2,-2) rectangle (2,2); \draw (-2,0) -- (2,0); \draw (0,-2)
-- (0,2); \draw[dotted,domain=-2:2] plot (\x, {\x});
\draw[domain=-2:0] plot (\x, {-4*(0.9*sqrt(abs(\x)/2)-0.5)});
\draw[domain=0:2] plot (\x, {4*(0.9*sqrt(abs(\x)/2)-0.5)});
\draw[dashed,domain=50/81:2, red] plot (\x,
{4*(0.9*sqrt(abs(4*(0.9*sqrt(abs(\x)/2)-0.5))/2)-0.5)});
\draw[dashed,domain=0:50/81, red] plot (\x,
{-4*(0.9*sqrt(abs(4*(0.9*sqrt(abs(\x)/2)-0.5))/2)-0.5)});
\draw[dashed,domain=-2:-50/81,red] plot (\x,
{-4*(0.9*sqrt(abs(-4*(0.9*sqrt(abs(\x)/2)-0.5))/2)-0.5)});
\draw[dashed,domain=-50/81:0,red] plot (\x,
{4*(0.9*sqrt(abs(-4*(0.9*sqrt(abs(\x)/2)-0.5))/2)-0.5)}); \draw
(-.396,-.396) rectangle (.396,.396);
\end{tikzpicture}
\caption{Lorenz map $f$ and its second iterate} \label{fig1}
\end{figure}

Let $f$ be an odd Lorenz map on $[-1,1]$. Observe that the second
iterate of Lorenz map, dashed graph in Figure \ref{fig1}, has two
fixed points $a$ and $-a$ and $f(a)=-a$. Choose $b\in(0,a)$ such
that $f^2(b)=-a$. We assume $f(1)>-f(b)$.
Now we follow the Bowen's construction of a fat horseshoe in \cite{B1}. For the sake of
completeness, we provide the details here. Let $(\beta_n)_{n\geq 0}$
be a sequence of positive numbers with
$$
\beta_0=2b,\quad\sum\limits_{n=0}^{\infty}\beta_n<2a\quad\text{and}\quad\beta_{n+1}/\beta_n\to
1.
$$
Denote a word of length $n$ with alphabet $\{0,1\}$ by
$w=w_1w_2\ldots w_n$ where $w_i\in\{0,1\}$ and $\ell(w)=n$. We are allowed to say that
the empty word $w=\emptyset$ has length 0 and
$I_{\emptyset}=[-a,a]$.  Let $I_{w0}$ and
$I_{w1}$ be the right and left intervals remaining after
removing the interior of $I_{w}^*$ from
$I_{w}$ where $I_{w}^*$ is the closed
interval of length $\beta_{\ell(w)}/2^{\ell(w)}$
and having the same center as $I_{w}$. Then
$$
K_I=\bigcap_{n=0}^{\infty}\bigcup_{\ell(w)=n}I_{w}
$$
is a Cantor set and
$\leb(K)=2a-\sum\limits_{n=0}^{\infty}\beta_n>0$. Now we replace the
function $f$ on $[f(b),-a]\cup[a,-f(b)]$ with new function and by
abuse of notation, we continue to write $f$ for the new function.
Let
$f:\bigcup_{w}f(I^*_{1w})\to\bigcup_{w}I^*_{w}$
be a function such that
\begin{itemize}
 \item $f|_{f(I^*_{1w})}$ is a $C^1$ orientation preserving diffeomorphism and $f(f(I^*_{1w}))=I^*_{w}$;
 \item $f'(f(x))f'(x)=2$ for any endpoint $x$ of $I^*_{1w}$;
 \item $\sup\limits_{x\in I^*_{1w}}|2-f'(f(x))f'(x)|\to 0$ as $\ell(w)\to\infty$.
\end{itemize}
Next we can extend $f$ continuously to $[f(b),-a]$ so that
$f:[f(b),-a]\to[-a,a]$ is a $C^1$ diffeomorphism with
$f'(f(x))f'(x)=2$ for $x\in K_I\cap[b,a]$. The map $f$ is defined on
$[a,-f(b)]$ by $f(x)=-f(-x)$. Therefore, if we assume that
$(f^2)'(a)=(f^2)'(b)=2$, then $f$ is a $C^1$-map such that
$f^2:[b,a]\to[-a,a]$ is a base map for a fat horseshoe in \cite{B1}.

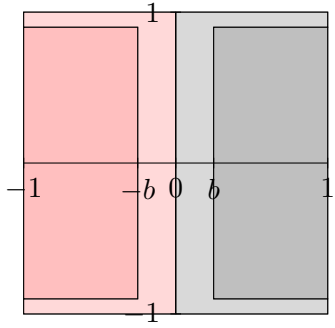
\begin{figure}[h]
\begin{center}
\begin{tikzpicture}[scale=2]

\draw (-1,-1) rectangle (1,1);

\filldraw [fill=gray!30] (0,-1) -- (0,1) -- (1,1) -- (1,.9) --
(.25,.9) -- (.25,-.9) -- (1,-.9) -- (1,-1) -- (0,-1);
\filldraw[fill=gray!50] (.25,-.9) rectangle (1,.9);

\filldraw [fill=red!15]
(0,1) -- (0,-1) -- (-1,-1) -- (-1,-.9) -- (-.25,-.9) --
(-.25,.9) -- (-1,.9) -- (-1,1) -- (0,1); \filldraw
[fill=red!25]
(-1,-.9) rectangle (-.25,.9);

\draw (-1,0) -- (1,0); \foreach \x/\xtext in { 0.25/b, 1,-0.25/-b,
-1,0} \draw (\x cm,1pt) -- (\x cm,-1pt) node[anchor=north]
{$\xtext$};

\draw (0,-1) -- (0,1); \foreach \y/\ytext in {-1,1} \draw (1pt,\y
cm) -- (-1pt,\y cm) node[anchor=east] {$\ytext$};

\end{tikzpicture}
\caption{Domain partition of the Poincar\'e map} \label{fig2}
\end{center}
\end{figure}

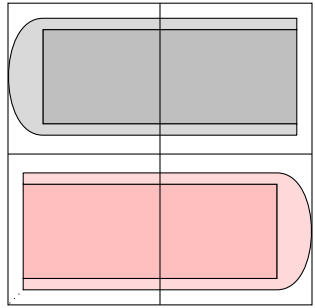
\begin{figure}[h]
\begin{center}
\begin{tikzpicture}[scale=2]
\draw[dotted,domain=-.99:-.77] plot (\x, {\x}); \draw (-1,-1)
rectangle (1,1);

\filldraw[fill=gray!50] (-.77,.2) rectangle (.9,.825);
\filldraw[fill=gray!30] (.9,.125) -- (-.77,.125) to[out=180,in=180]
(-.77,.9) -- (.9,.9) -- (.9, .825) -- (-.77,.825) -- (-.77,.2) --
(.9,.2) -- (.9,.125);

\filldraw[fill=red!25]
(-.9,-.825) rectangle (.77,-.2);
\filldraw[fill=red!15]
(-.9,-.125) -- (.77,-.125) to[out=0,in=0] (.77,-.9) --
(-.9,-.9) -- (-.9, -.825) -- (.77,-.825) -- (.77,-.2) -- (-.9,-.2)
-- (-.9,-.125);

\draw (-1,0) -- (1,0);

\draw (0,-1) -- (0,1);

\end{tikzpicture}
\caption{Image of the Poincar\'e map} \label{fig3}
\end{center}
\end{figure}

We can now define a Poincar\'e map $F(x,y)=(f(x),g(x,y))$ on the
section $[-1,1]\times [-1,1]$. Notice that $F$ is not defined on
line $0\times [-1,1]$. The domain of $F$ is partitioned like Figure
\ref{fig2}, and $g$ is a $C^1$-map such that
$$
g(x,y)=\sign(x)f^{-1}(\sign(x)y)
$$
for $(x,y)\in[-1,-b]\cup[b,1]\times[f(b)-\epsilon,-f(b)+\epsilon]$,
where $\epsilon=\frac{f(1)+f(b)}{2}$, $\sign(x)=|x|/x$ and $f^{-1}$ is the inverse of
right branch of $f$.
In rest of the domain $g$ is defined in a way
that the image of these areas under $F$ is shown in Figure
\ref{fig3}.
For $(x,y)\in [-a,-b]\cup[b,a]\times[-a,a]$, $F^2(x,y)$ is given by
$$
F^2(x,y)=\left(f^2(x),-\sign(x)f^{-1}(-f^{-1}(\sign(x)y))\right).
$$
Let $A$ be the square $[-a,a]\times[-a,a]$, then
$$
H=\bigcap_{k=-\infty}^{\infty}F^{2k}(A)=K_I\times K_I,
$$
and $\leb_2(H)=\leb(K_I)^2>0$. On the other hand,
$$
H\subseteq \bigcap_{k=0}^{\infty}F^{2k}(A)\subseteq
\bigcap_{k=0}^{\infty}F^{k}(\Sigma).
$$
The last inclusion holds since $f^{2k}(\Sigma)\subseteq
f^{2k-1}(\Sigma)$.
\section{Appendix}
Our method provides another proof for the fact that a $C^1$-generic diffeomorphism does't admit partially hyperbolic horseshoe-like attractors of positive measure. We need some other notions. An invariant
set $\Lambda$ admits a dominated splitting $E\oplus F$ over
$\Lambda$ if there is $0<\lambda<1$ such that
$$\|Df|_{E(x)}\|~\|Df^{-1}|_{F(f(x))}\|\leq \lambda.$$ A dominated
splitting $E^s\oplus F$ is partially hyperbolic if
$\|Df|_{E^s}\|<\lambda$. It is known that if $f$ has a partially
hyperbolic set $\Lambda_f$ then there is a neighborhood $U$ of
$\Lambda_f$ and $\mathcal{U}$ of $f$ in $C^r (r\geq 1)$ topology such that any
$g\in\mathcal{U}$ admits a partially hyperbolic structure on the
invariant set $\Lambda_g=\cap_{n\in\mathbb{Z}} g^n(U)$, see \cite{hps}.
The set $\Lambda_g$ is not necessarily close to the initial set
$\Lambda_f$ in the Hausdorff metric. As in the case of generalized
hyperbolic attractors, there is a local dynamical lamination
$W^s_\delta$ of non-uniform size which integrates the direction $E^s$,
that is $T_x W^s_\delta(x)=E^s(x)$. Furthermore, for any
sufficiently $C^1$-close diffeomorphism $g$ and any two points
$x\in\Lambda_f$ and $y\in\Lambda_g$, the two local stable manifolds
$W^s_\delta(x,f)$ and $W^s_\delta(x,g)$ are sufficiently close in
Hausdorff metric.
As before, we say that a partially hyperbolic set $\Lambda$ is {\it
horseshoe-like} if it does not contain any local stable manifold.
\begin{theorem*}
\label{thma} There is a residual subset $\mathcal{R}_0\subset \diff^1(M)$ such that for any $f\in\mathcal{R}_0$, if $\Lambda$ is a
partially hyperbolic horseshoe-like attractor of $f$ then
$Leb(\Lambda)=0$.
\end{theorem*}
\begin{proof}
Let $\mathcal{B}$ be a countable open base for the topology of $M$ and
$\{V_n\}_{n\in\mathbb{N}}$ be the family of all finite union of the
elements in $\mathcal{B}$. Suppose that $\mathcal{V}_n$ is the set of all
$f\in \diff^1(M)$ such that $f(\overline{V_n})\subset V_n$ and
$\cap_{m\in\mathbb{N}} f^m(V_n)$ is partially hyperbolic. By the
robustness of the partial hyperbolicity, the $\mathcal{V}_n$ is an
open set. Now, define $\mathcal{F}_n:\mathcal{V}_n\to 2^M$ by 
$\mathcal{F}_n(f)=\cap_{m\in\mathbb{N}} f^m(V_n)$. It is not
difficult to see that $\mathcal{F}_n$ is upper semi-continuous.
Hence, there is a residual subset $\mathcal{R}_n$ of $\mathcal{V}_n$
such that $\mathcal{F}_n$ varies continuously on it. We note that by
the regularity of Lebesgue measure, any $f\in \mathcal{R}_n$ is also
a continuity point of the mapping $\mathcal{G}_n:\mathcal{V}_n\to
\mathbb{R}$ given by $\mathcal{G}_n(f)=Leb (\mathcal{F}_n(f))$. Now, for any $n$, put
$\mathcal{H}_n=\diff^1(M)\setminus \overline{\mathcal{V}_n}$,
$\mathcal{R}_n^\prime=\mathcal{H}_n\cup \mathcal{R}_n$ and $\mathcal{R}=\cap_{n\in\mathbb{N}} \mathcal{R}_n^\prime$. Note that $\mathcal{R}_n^\prime$ is a residual subset of $\diff^1(M)$. Let $f\in\mathcal{R}$ has a partially hyperbolic horseshoe-like attractor $A=\cap_{n\in\mathbb{N}} f^n(U)$. Choose $n_0$ such that
$A\subset V_{n_0}\subset U$. It can be easily seen that $f\in
\mathcal{V}_{n_0}$ and $A=\mathcal{F}_{n_0}(f)$. Since $f\in \mathcal{R}$, we have $f\in
\mathcal{R}_{n_0}^\prime$ 
and thus
$f\in\mathcal{R}_n$. Now, let $(f_n)_n$ be a sequence of
$C^2$-diffeomorphisms such that $f_n\to f$. Since $\mathcal{F}_{n_0}(f)$ is
horseshoe-like, the same holds for $\mathcal{F}_{n_0}(f_n)$, for sufficiently
large $n$. To show it, let assume by contradiction that $\mathcal{F}_{n_0}(f_{n_k})$ contains a local stable
manifold $W^s_\delta(x_{n_k})$ for a subsequence $(n_k)_k$. Without loss of generality, we may assume that the subsequence $(n_k)_k$ converges to $x$. 
By the continuity of the stable manifolds, we
have
$$W^s_\delta(x)=\lim_{k\to\infty}W^s_\delta(x_{n_k},f_{n_k})\subset
\lim_{k\to\infty}\mathcal{F}_{n_0}(f_{n_k})= \mathcal{F}_{n_0}(f),$$ and hence $W^s_\delta(x)\subset
\mathcal{F}_{n_0}(f)$ which is a contradiction. By theorem above,
$Leb(\mathcal{F}_{n_0}(f_{n_k}))=0$ and hence,
\[Leb(A)=\mathcal{G}_{n_0}(f)=\lim_{k\to\infty} \mathcal{G}_{n_0}(f_{n_k})=0.\]
\end{proof}
\section*{Acknowledgments}
During the preparation of this article the first author was partially
supported by grant from IPM (No. 96370119). He also thanks ICTP for supporting through the association schedule.
{}
\end{document}